\documentclass[12pt]{article}
\usepackage{amssymb}
\usepackage{array}
\usepackage{latexsym}
\usepackage{enumerate}
\usepackage{amsmath}
\usepackage{amsfonts}
\usepackage{amsthm}
\usepackage{graphicx}
\usepackage{psfrag}
\usepackage{comment}
\usepackage[english]{babel}
\usepackage[T1]{fontenc}
\usepackage{geometry}
\usepackage{hyperref}
\ifx\smallsetminus\undefined
\def\smallsetminus{\setminus}
\fi
\title{An alternative construction of B-M and B-T unitals in desarguesian planes}
\author{A.~Aguglia\footnote{Dipartimento di Matematica, Politecnico di Bari,
  Via Amendola 126/B, I-70126 Bari}
,\,  L.~Giuzzi\footnote{Dipartimento di Matematica, Facolt\`a di Ingegneria,
  Universit\`a degli Studi di Brescia, Via Valotti 9, I-25133
  Brescia}, \,G.~Korchm\'aros\footnote{Dipartimento di Matematica
  e Informatica,
 Universit\`a degli Studi della Basilicata, Campus Macchia Romana,
Viale dell'Ateneo Lucano, 10, I-85100 Potenza}}
\date{}
\theoremstyle{plain}
\newtheorem{prop}{Proposition}[section]
\newtheorem{theorem}[prop]{Theorem}

\newtheorem{lemma}[prop]{Lemma}
\theoremstyle{definition}
\newtheorem{remark}[prop]{Remark}

\newcommand{\cD}{\mathcal D}

\def\cU{\mathcal U}

\def\cP{\mathcal P}
\def\cL{\mathcal L}

\def\cH{\mathcal H}

\newcommand{\cC}{\mathcal C}

\newcommand{\fA}{\mathfrak{A}}
\newcommand{\fC}{\mathfrak{C}}
\newcommand{\PG}{\mathrm{PG}}
\newcommand{\AG}{\mathrm{AG}}
\newcommand{\GF}{\mathrm{GF}}

\newcommand{\Tr}{\mathrm{Tr}\,}
\begin{document}

\maketitle

\begin{abstract}
We present a new construction of non-classical unitals from a
classical unital $\cU$ in $\PG(2,q^2)$. The resulting
non-classical unitals are B-M unitals. The idea is to find a
non-standard model $\Pi$ of $\PG(2,q^2)$ with the following three
properties:
\begin{enumerate}[(i)]
\item points of $\Pi$ are those of $\PG(2,q^2)$;
\item
  lines of $\Pi$ are certain lines and conics of $\PG(2,q^2)$;
\item
  the points in $\cU$ form a non-classical B-M unital in $\Pi$.
\end{enumerate}
Our construction also works for the B-T unital, provided that conics are
replaced by certain algebraic curves of higher degree.
\end{abstract}


\section{Introduction}
A classical unital $\cU$ in the Desarguesian plane $\PG(2,q^2)$ is
the set of all absolute points of a non-degenerate unitary
polarity. Up to a projectivity of $\PG(2,q^2)$, $\cU$ consists of
all the $q^3+1$ points of the non-degenerate Hermitian curve $\cH$
with equation $y^q+y-x^{q+1}=0$. The relevant combinatorial
property of $\cU$, leading to important applications in coding
theory, is that $\cU$ is a {\em two-character set} with parameters
$1$ and $q+1$, that is, a line in $\PG(2,q^2)$ meets $\cU$ in
either $1$ or $q+1$ points. A {\em unital} in $\PG(2,q^2)$ is
defined by this combinatorial property, namely it is a
two-character set of size $q^3+1$ with parameters $1$ and $q+1$.

The known non-classical unitals are the B-M unitals due to
Buekenhout and Metz, see \cite{B,M}, and the B-T unitals due to
Buekenhout, see \cite{B}. They were constructed by an ingenious
idea, relying on the Bruck--Bose representation of $\PG(2,q^2)$ in
$\PG(4,q)$ and exploiting properties of spreads and ovoids (in
particular, quadrics). For $q$ odd, an alternative construction
for special B-M unitals which are the union of $q$ conics sharing
a point has been given by Hirschfeld and Sz\H{o}nyi \cite{hsz}.
Such B-M unitals are called H-Sz type B-M unitals.

In this paper, we present a new construction for non-H-Sz type B-M
unitals. The key idea, as described in Abstract, is fully realised
within $\PG(2,q^2)$, and it uses only quadratic transformations.
This method also works for B-T unitals, provided that quadratic
transformations are replaced by certain birational transformations.

Our notation and terminology are standard. For generalities on
unitals in projective planes the reader is referred to
\cite{GLEB,EG1,EG}. Basic facts on rational transformations of
projective planes are found in \cite[Section 3.3]{HKT}.

\section{A non-standard model of $\PG(2,q^2)$}
\label{nonst} Fix a projective frame in $\PG(2,q^2)$ with
homogeneous coordinates $(x_0,x_1,x_2)$, and consider the affine
plane $\AG(2,q^2)$ whose infinite line $\ell_{\infty}$ has
equation $x_0=0$. Then $\AG(2,q^2)$ has affine coordinates $(x,y)$
where $x=x_1/x_0,\,y=x_2/x_0$ so that $X_\infty=(0,1,0)$ and
$Y_\infty=(0,0,1)$ are the infinite points of the horizontal and
vertical lines, respectively.

Fix a non-zero element $a\in\GF(q^2)$. For $m,d\in\GF(q^2)$ and
$a\in\GF(q^2)^*$, let $\cC_a(m,d)$ denote the parabola of equation
$y=ax^2+mx+d$. Consider the incidence structure $\fA_a=(\cP,\cL)$
whose points are the points of $\AG(2,q^2)$ and whose lines are
the vertical lines of equation $x=k$, together with the parabolas
$\cC_a(m,d)$ where $m,d,k$ range over $\GF(q^2)$.
\begin{lemma}
\label{lemm1} For every non-zero $a\in \GF(q^2)$, the incidence
structure $\fA_a=(\cP,\cL)$ is an affine plane isomorphic to
$\AG(2,q^2)$.
\end{lemma}
\begin{proof}
The birational transformation $\varphi$ given by
 \begin{equation} \label{fi}\varphi: (x,y)\mapsto (x,y-ax^2), \end{equation}
transforms vertical lines  into themselves, whereas the generic
line $y=mx+d$ is mapped into the parabola $\cC_a(m,d)$. Therefore,
$\varphi$ determines an isomorphism
\[\fA_a \simeq \AG(2,q^2),\] and the assertion is proven.
\end{proof}

Completing $\fA_a$ with its points at infinity in the usual way
gives a projective plane isomorphic to $\PG(2,q^2)$. Note that the
infinite point $Y_\infty$ of the vertical lines of $\AG(2,q^2)$ is
also the infinite point of the vertical lines of $\fA_a$.

For $q$ an odd power  of $2$, a different, yet similar,
construction will also be useful in our investigation. The
construction depends on some known facts about Galois fields of
even characteristic. Let $\varepsilon \in \GF(q^2) \setminus
\GF(q)$ such that $\varepsilon^2+\varepsilon+\delta=0$, for some
$\delta \in \GF(q)\setminus \{1\}$ with $\Tr(\delta)=1$. Here, as
usual, $\Tr$ stands for the trace function $\GF(q)\to \GF(2)$.
Then $\varepsilon^{2q}+\varepsilon^q+\delta=0$. Therefore,
$(\varepsilon^q+\varepsilon)^2+(\varepsilon^q+\varepsilon)=0$,
whence $\varepsilon^q+\varepsilon+1=0$. Moreover, if $q$ is an odd
power of $2$, then
$$\sigma: x \mapsto x^{2^{(e+1)/2}}$$
is an automorphism of $\GF(q)$.

For any $m,d\in\GF(q^2)$ let $\cD(m,d)$ denote the plane algebraic
curve of equation
\begin{equation}
\label{gam} y=[((x^q+x)\varepsilon+x)^{\sigma+2}+(x^q+x)^\sigma+
((x^q+x)\varepsilon+x)(x^q+x)]\varepsilon+bx^{q+1} +mx+d
\end{equation}
where $b$ is a given element in $\GF(q^2)\setminus \GF(q)$.

Introduce the incidence structure $\fA_\varepsilon'=(\cP',\cL')$
whose points are the points of $\AG(2,q^2)$ and whose lines are
the vertical lines of equation $x=k$, together with the curves
$\cD(m,d)$ where $m,d,k$ range over $\GF(q^2)$.
\begin{lemma}
\label{lemm2} The incidence structure
$\fA_\varepsilon'=(\cP',\cL')$ is an affine plane isomorphic to
$\AG(2,q^2)$.
\end{lemma}
\begin{proof}
The argument  in the
proof of Lemma \ref{lemm1} works also in this case,
provided that $\varphi$ is replaced
by the birational transform $\gamma$ defined by
 $$\gamma:(x,y)\mapsto
 (x,y+[((x^q+x)\varepsilon+x)^{\sigma+2}+(x^q+x)^\sigma+
((x^q+x)\varepsilon+x)(x^q+x)]\varepsilon+bx^{q+1})$$
 \end{proof}
\section{The Construction}
Before presenting our construction we recall the equations of B-M
unitals and B-T unitals in $\PG(2,q^2)$.
\begin{prop}{\rm(Baker and Ebert\cite{BE2}, Ebert \cite{BE3,EG1}).}
\label{BM} For $a,b\in\GF(q^2)$, the point--set
\[U_{a,b}=\{(1,x,ax^2+bx^{q+1}+r)|x\in\GF(q^2), r \in\GF(q)\}\cup \{Y_\infty\}\]
 is a B-M unital in $\PG(2,q^2)$ if and only if Ebert's discriminant
 condition is satisfied, that is
 for odd $q$,
\begin{itemize}
\item[\rm(i)]
 $4a^{q+1}+(b^q-b)^2$ is a non--square in $\GF(q)$,
\end{itemize}
 and for $q$ even,
 \begin{itemize}
\item[\rm(ii)]
 $b\notin \GF(q)$ and $ \Tr(a^{q+1}/(b^q+b)^2)=0$.
\end{itemize}
Conversely, every B-M unital has a representation as $U_{a,b}$.
\end{prop}
\begin{prop}
\label{ebertbis} With the above notation,
    \begin{enumerate}
    \item[\rm(i)] $U_{a,b}$ is classical if and only if $a=0$;
    \item[\rm(ii)] $U_{a,b}$ is a H-Sz type B-M unital if
    and only $a^{(q+1)/2}\in \GF(q^2)\setminus\GF(q)$ and $b\in \GF(q)$.
    \end{enumerate}
\end{prop}
\begin{proof} This is a direct corollary of \cite[Theorems 1 and
12]{EG1}.
\end{proof}
\begin{prop}{}
\label{xEp}
Let $q=2^e$, where $e>1$ is an odd integer. In the
above notation, the point--set
\begin{equation}
\label{nuovo3}
\begin{array}{l}
U_\varepsilon=\{(1,x,[((x^q+x)\varepsilon+x)^{\sigma+2}+(x^q+x)^\sigma+
((x^q+x)\varepsilon+x)(x^q+x)]\varepsilon+r \\ \mid x\in
\GF(q^2),r\in\GF(q) \} \cup \{Y_\infty\},
\end{array}
\end{equation}
is a  B-T unital in $\PG(2,q^2)$. Conversely, every B-T unital may be
represented as $U_{\varepsilon}$ for some choice of $\varepsilon$.
\end{prop}
\begin{proof}
From \cite{EGL,EG1}, the point--set
 \begin{equation}
 \label{nuovo4}
 U_\varepsilon=\{(1,s+t \varepsilon, (s^{\sigma+2}+t^{\sigma}+st)
 \varepsilon+r |r,s,t \in\GF(q) \} \cup \{Y_\infty\}
\end{equation}
is a B-T unital and, conversely, every B-T unital has such an
equation. Let $x=s+t\varepsilon$. Then, $t=x^q+x$ and
$s=x+(x^q+x)\varepsilon$. Substituting $t$ and $s$ in
\eqref{nuovo4} gives the result.
\end{proof}
If $b\in \GF(q^2)\setminus\GF(q)$ then, from
Proposition \ref{ebertbis}, the point--set
\begin{equation}
\label{nuovo2} \cU_b=\{(1,x,bx^{q+1}+r)|x\in\GF(q^2), r \in
\GF(q)\}\cup \{Y_\infty\}
\end{equation}
is a classical unital in $\PG(2,q^2)$. As pointed out in
Section \ref{nonst}, $\cU_b$ can be regarded as a point--set in the
projective closure of $\fA_a$ and, for $q$ even, also as a
point--set of the projective closure of $\fA_{\varepsilon}'$. The question
arises whether $\cU_b$ is still a unital in these planes. Our main
result, stated in the following two theorems, shows that the
answer is positive.
\begin{theorem}\label{teo} Let $a\in \GF(q^2),\,b\in
\GF(q^2)\setminus \GF(q)$. If $(a,b)$ satisfies Ebert's
discriminant condition, then $\cU_b$ is the non--classical B-M
unital $U_{-a,b}$ in the projective closure of $\fA_a$.
Conversely, every non-H-Sz type B-M unital  is obtained in
this way.
\end{theorem}
\begin{proof}
Let $P=(\xi,\eta)$ an affine point in $\fA_a$. This point, viewed as
an element of $\AG(2,q^2)$, has coordinates $x=\xi$ and
$y=\eta+a\xi^2$. From \eqref{nuovo2},
\begin{equation}
\label{rqnuovo1} \cU_b=\{(1,\xi,-a\xi^2+ b\xi^{q+1}+r\mid \xi\in
\GF(q^2), r \in\GF(q)\}\cup \{Y_\infty\}.
\end{equation}
This shows that $\cU_b$ and $U_{-a,b}$ coincide in $\fA_a$. Since
$(-a,b)$ also satisfies Ebert's discriminant condition, $U_{-a,b}$
is a B-M unital in the projective closure of $\fA_a$. By
Proposition \ref{ebertbis}, $U_{-a,b}$ is a non-H-Sz type B-M unital.
\end{proof}
\begin{theorem}\label{teo1}
Let $q=2^e$, with $e>1$ an odd integer.  Then $\cU_{\varepsilon}$
is a non--classical B-T unital in the projective closure of
$\fA_{\varepsilon}'$.
\end{theorem}
\begin{proof}
We use the same argument as in the preceding proof. The point
$P=(\xi,\eta)$ of $\fA_{\varepsilon}'$, viewed as an element of $\AG(2,q^2)$, has
coordinates $x=\xi$ and
$$y=\eta+[((\xi^q+\xi)\varepsilon+\xi)^{\sigma+2}+(\xi^q+\xi)^\sigma+
((\xi^q+\xi)\varepsilon+\xi)(\xi^q+\xi)]\varepsilon+b\xi^{q+1}.$$
From \eqref{nuovo2},
\begin{equation*}
\begin{array}{l}
\cU_b=\{(1,\xi,[((\xi^q+\xi)\varepsilon+x)^{\sigma+2}+(\xi^q+\xi)^\sigma+
((\xi^q+\xi)\varepsilon+\xi)(\xi^q+\xi)]\varepsilon+r \\ \mid
\xi\in\GF(q^2),r\in\GF(q)\} \cup \{Y_\infty\}.
\end{array}
\end{equation*}
By Proposition \ref{xEp} we have that $\cU_b$ and $\cU_{\varepsilon}$ coincide in $\fA_{\varepsilon}'$ and the assertion follows.
\end{proof}

\subsection{ An alternative proof of Theorem \ref{teo}}

The above proofs of Theorem \ref{teo} and  \ref{teo1} depend on
the explicit equations for B-M and B-T unitals, as given in
Propositions \ref{BM} and \ref{xEp}. Here we provide a direct
proof of Theorem \ref{teo}. Without loss of generality, we assume
that $q\geq 3$.

Let $\cH$ be the set of all points in $\AG(2,q^2)$ of the affine
Hermitian curve $\cC$ of equation
\begin{equation}
\label{nuovo6} y^q-y+(b-b^q)x^{q+1}=0,\, \,b\not\in \GF(q),
\end{equation}
Then, $\cH\cup \{Y_\infty\}$ is a classical unital in $\PG(2,q^2)$.
We prove that $\cH\cup \{Y_\infty\}$ is also a unital in the
projective closure of $\fA_a$.

We first need the following lemma.
\begin{lemma}\label{lem1}
For every $m,d\in \GF(q^2)$, the parabola $\cC_a(m,d)$ and $\cH$
have either $1$ or $q+1$ points in  $\AG(2,q^2)$.
\end{lemma}
\begin{proof}
The number of solutions
  $(x,y)\in\GF(q^2)\times \GF(q^2)$ of
  the system
  \begin{equation}
    \left\{\begin{array}{l}\label{sis1}
        y^q-y+(b-b^q)x^{q+1}=0 \\
        y-ax^2-mx-d=0
      \end{array}\right.
  \end{equation}
gives the number of points in common of $\cH$ and $\cC_a(m,d)$. To
solve this system, recover the value of $y$ from the second
equation and substitute it in the first. The result is
\begin{equation}
\label{eq1}
a^qx^{2q}+(b-b^q)x^{q+1}+m^qx^{q}-ax^2-mx+d^q-d=0.
\end{equation}
Consider now $\GF(q^2)$ as a vector space over $\GF(q)$, fix a
basis $\{1,\varepsilon\}$ with $\varepsilon\in\GF(q^2)\setminus
\GF(q)$, and write the elements in $\GF(q^2)$ as a linear
combination with respect to this basis, that is,
$z=z_0+z_1\varepsilon$, with $z\in \GF(q^2)$ and
$z_0,z_1\in\GF(q)$. Thus, \eqref{eq1} becomes an equation over
$\GF(q)$. We investigate separately the even $q$ and odd $q$
cases.

For $q$ even, $\varepsilon$ may be chosen  as in Section \ref{nonst}. With
this choice of $\varepsilon$, \eqref{eq1} becomes
 \begin{equation}\label{eq2}
(a_1+b_1)x_0^2+[(a_0+a_1)+\nu(a_1+b_1)]x_1^2+b_1x_0x_1+m_1x_0+(m_0+m_1)x_1+d_1=0.
 \end{equation}
We shall represent the
the solutions $(x_0,x_1)$ of \eqref{eq2} as points of the affine plane
$\AG(2,q)$ over $\GF(q)$ arising from the vector space $\GF(q^2)$.
In fact, \eqref{eq2} turns out to be the equation of
a (possibly degenerate) affine conic $\Xi$ of $\AG(2,q)$.
Actually, $\Xi$ is either an ellipse or is a single point. To
prove this, we have to show that it has no point  at infinity;
that is, we need to prove that the points $P=(x_0,x_1,0)$ with
\begin{equation}\label{eq3}
(a_1+b_1)x_0^2+[(a_0+a_1)+\nu(a_1+b_1)]x_1^2+b_1x_0x_1=0,
\end{equation}
do not lie in $\PG(2,q)$. Obviously, this is the case if and only
if \eqref{eq3} admits only the trivial solution over $\GF(q)$. A
necessary a sufficient condition for this is
\begin{equation}
\label{nuovo11} \Tr\left(\frac{(a_1+b_1)
[(a_0+a_1)+\nu(a_1+b_1)]}{b_1^2}\right)=1.
\end{equation}
In our case, \eqref{nuovo11} holds as it
follows directly from Ebert's discriminant condition; see \cite[page 83]{GLEB}.  Therefore, $\Xi$ is
either an ellipse or it consists of a single point; hence,
$\cC_a(m,d)$ meets $\cH$ in either $1$ or $q+1$ points.

For $q$ odd, an analogous argument is used. For this purpose, as
in \cite{EG1}, choose a primitive element $\beta$ of $\GF(q^2)$
and let $\varepsilon=\beta^{(q+1)/2}$. Then,
$\varepsilon^q=-\varepsilon$ and $\varepsilon^2$ is a primitive
element of $\GF(q)$. With this choice of $\varepsilon$,
\eqref{eq1} becomes
\begin{equation}
\label{eqodd1} (b_1+a_1)\varepsilon^2
x_1^2+2a_0x_0x_1+(a_1-b_1)x_0^2+m_0x_1+m_1x_0+d_1=0.
\end{equation}
The discussion of the (possibly degenerate) affine conic $\Xi$ of
equation \eqref{eqodd1} may be carried out exactly as in the even order
case. The points $P=(x_0,x_1,0)$ of $\Xi$ at infinity are
determined by
\[ (b_1+a_1)\varepsilon^2 x_1^2+2a_0x_0x_1+(a_1-b_1)x_0^2=0, \]
and this equation has only the trivial solution over $\GF(q)$,
since Ebert's discriminant condition implies that $4a^q+(b^q-b)^2$
is non-square in $\GF(q)$.
\end{proof}
Lemma \ref{lem1} together with \cite[Theorem 12.16]{H} have the
following corollary.
\begin{theorem}
\label{nuovo7} The point-set $\cH\cup \{Y_\infty\}$ is a unital in
the projective closure of $\fA_a$.
\end{theorem}
To show that $\cH\cup \{Y_\infty\}$ is a non-classical unital in the
projective closure of $\fA$, we rely on some elementary facts on
algebraic curves.
\begin{lemma}
\label{nuovo30} The points of $\cH$ in $\fA_a$ lie on the
absolutely irreducible affine plane curve $\cC'$ of equation
$$\eta^q-\eta+(b-b^q)\xi+a^q\xi^{2q}-a\xi^2=0.$$
\end{lemma}
\begin{proof}
  The plane curve $\cC'$ is absolutely irreducible, see \cite[Lemma 12.1]{HKT}.
  If $P=(\xi,\eta)$ is a point of $\cH$ in $\fA_a$, then $P$, regarded
as a point of $\AG(2,q^2)$, has coordinates $x,y$ with
$x=\xi,\,y=\eta+a\xi^2$. Since $(x,y)$ satisfies  \eqref{nuovo6},
$$(\eta+a\xi^2)^q-\eta-a\xi^2+(b-b^q)\xi=0$$ holds.
This implies that $P=(\xi,\eta)$ is a point of $\cC'$.
\end{proof}
\begin{theorem}
\label{nuovoth1} The point-set $\cH\cup \{Y_\infty\}$ is a
non-classical unital in the projective closure of $\fA_a$.
\end{theorem}
\begin{proof} Assume, on the contrary,
  that $\cH$ coincides in $\fA_a$ with the point--set of
a non-degenerate affine Hermitian curve $\cD'$. Then, $\cC'$ and
$\cD'$ have at least $q^3$ common points. Since $\deg \cC'=2q$ and
$\deg \cD' =q+1$ and $2q(q+1)<q^3$, B\'ezout's theorem, see
\cite[Theorem 3.13]{HKT}, implies that $\cC'$ and $\cD'$ share a
common component. This contradicts Lemma  \ref{nuovo30}.
\end{proof}
Finally, we prove that $\cH\cup \{Y_\infty\}$ is a B-M unital in
the projective closure of $\fA_a$. Our proof relies on the
Ebert-Wantz group-theoretic characterization of B-M unitals of a
Desarguesian plane: A unital $\cU$ of $\PG(2,q^2)$ is a B-M unital
if, and only if, $\cU$ is preserved by a linear collineation group
of order $q^3(q-1)$ which is the semidirect product of a subgroup
$S$ of order $q^3$ by a subgroup $R$ of order $q-1$. Moreover, $S$
is Abelian if, and only if, $\cU$ is a H-Sz type B-M unital; see
\cite{EW} and \cite[Theorem 12]{EG1}. For more results on the
collineation group of a B-M unital, see \cite{AL1,AL2}.
\begin{theorem}
\label{thnuovo3} In the projective closure of $\fA_a$, the
point-set $\cH\cup \{Y_\infty\}$ is a non-Sz-H type B-M unital.
\end{theorem}
\begin{proof}
A straightforward computation shows that for any point $P=(u,v)\in
\cH$ in $\fA_a$ and for any $\lambda\in \GF(q)^*$, the affinities
\begin{equation}
\label{nuovo33}
\begin{array}{llll}
\alpha_{u,v}:&\, (\xi,\eta)&\to&(\xi+u,\eta-2au\xi+u^q(b-b^q)\xi+v),\\
\beta_{\lambda}:&\, (\xi,\eta)&\to&(\lambda \xi, \lambda^2\eta)
\end{array}
\end{equation}
of $\fA_a$ preserve $\cH$. The group $S$ of the linear
collineations $\alpha_{u,v}$ with $P=(u,v)$ ranging over $\cH$ is a
non-Abelian group of order $q^3$. Write $R$ for the group of the linear
collineations $\beta_{\lambda}$ as $\lambda$  varies $\GF(q)^*$.
It turns out that the group $G$ generated by all these collineations
has order
$q^3(q-1)$ and is the semidirect product $S\rtimes R$, and the assertion follows from the
Ebert-Wantz characterization.
\end{proof}
\begin{remark} Theorem \ref{nuovoth1} may also be
proven without using algebraic geometry. The idea is to write the
equation of the tangent parabolas $\cC_a(m,d)$ at the points of
the classical unital $\cH\cup\{Y_{\infty}\}$ and use Thas'
characterization \cite{thas} involving the feet of a point on a
unital. If $P=(w,z)\in\cH$ then the unique tangent parabola to
$\cH$ at $P$ has equation
\begin{equation}
\label{tanp} y=ax^2+\left(-2aw+(b-b^q)w^q\right)x-z^q+aw^2.
\end{equation}
For $q$ odd, Theorem \ref{thnuovo3} can also be proven replacing group theoretic
arguments with some geometric characterisations results
depending on Baer sublines, due to Casse, O'Keefe, Penttila and Quinn;
see \cite{casse,qc} and \cite[Theorem 11]{EG1}.
\end{remark}
\section{Absolutely
  irreducible curves containing all points of a B-M unital in $\PG(2,q^2)$}
For $a,b\in \GF(q^2)$ satisfying Ebert's discriminant condition,
the absolutely irreducible plane curve $\Gamma_{a,b}$ of
$\PG(2,q^2)$ with affine equation
\begin{equation}
\label{tit4} y^q-y-a^qx^{2q}+ax^{2}+(b-b^q)x^{q+1}=0.
\end{equation}
contains all points of the unital $U_{-a,b}$. We prove some
properties of $\Gamma_{a,b}$.
\begin{theorem}\label{teo3}
The curve $\Gamma_{a,b}$ is birationally equivalent over $GF(q^2)$
to a non--degenerate Hermitian curve.
\end{theorem}
\begin{proof} The birational map $(x,y)\to (x,y-ax^2)$ transforms $\Gamma_{a,b}$ into the Hermitian curve
$\cC$ of equation \eqref{nuovo6}.
\end{proof}
\begin{theorem}
\label{thnuovo4} $\Gamma_{a,b}$ is the unique plane curve of
minimum degree which contains all the points of the B-M unital
$U_{-a,b}$ in $\PG(2,q^2)$.
\end{theorem}
\begin{proof}
Let $\Psi$ be (a not necessarily absolutely irreducible) plane
curve of $\PG(2,q^2)$ of degree $d\leq 2q$ which contains all
points of $U_{-a,b}$. Obviously, $\Gamma_{a,b}$ and $\Psi$ have at
least $q^3+1$ common points. From B\'ezout's theorem \cite[Theorem
3.13]{HKT}, $\Gamma_{a,b}$ is a component of $\Psi$. Since $\deg\,
\Gamma_{a,b}\geq \deg\, \Psi$, this is only possible when they
coincide.
\end{proof}

\begin{remark}
In 1982,  Goppa introduced a general construction technique for
linear codes from algebraic curves defined over a finite field;
see \cite{Gp}.
In the current literature, these codes are called {\em algebraic-geometry}.

The parameters of linear codes arising from  a Hermitian curve by
Goppa's method  were computed in \cite{vL}. These codes turn out
to perform very well, when compared with  Reed-Solomon codes of
similar length and dimension.

In \cite{EG1}, Ebert raised  the question whether the parameters
of the  codes arising from $\Gamma_{a,b}$   by  Goppa's
construction were close to maximum distance separable codes.

Since the algebraic-geometric codes are determined by the
function fields of the related algebraic curves and the function
fields of  two birationally equivalent plane curves are
isomorphic, Theorem \ref{teo3} implies that the algebraic-geometry
codes arising from  the Hermitian curve $\cC$ and those arising
from the  curve $\Gamma_{a,b}$ are the same.
\end{remark}

\section{B-M unitals and cones of $\PG(3,q^2)$}

We  present another way to  construct a non-classical B-M unital
using a Hermitian curve and a suitable cone of $\PG(3,q^2)$.

Let $x_0, x_1, x_2, x_3$ denote homogeneous coordinates in
$\PG(3,q^2)$. Consider the Hermitian curve
$\cH=\{(1,t,bt^{q+1}+r)|t\in\GF(q^2), r \in
\GF(q)\}\cup\{Y_{\infty}\} $ and the map $\phi: \cH \mapsto
\PG(3,q^2)$ which transforms the point $P(1,t,bt^{q+1}+r)$ into the point
$\phi(P)=(1,t, t^2, bt^{q+1}+r)$ and $Y_{\infty}=(0,0,1)$  into
$\phi(Y_{\infty})=(0,0,0,1)$.

The map $\phi$ is injective; thus, the set $\phi(\cH)$ consists of
$q^3+1$ points lying on the cone $\fC$ represented by
$x_0x_2=x_1^2$. The point $Q(0,0,1,-a)$ does not lie on the cone
$\fC$; hence, the projection $\rho$  from $Q$ to the plane $\pi:
x_2=0$ is well defined. The point  $\phi(Y_{\infty})$ is on $\pi$
thus  we get $\rho(0,0,0,1)=(0,0,0,1)$.

For any $(t,r) \in\GF(q^2) \times\GF(q)$, set
$P_{t,r}=(1,t,bt^{q+1}+r)$. The line $P_{t,r}Q$ has point set
\[\{(1,t,t^2+\lambda, bt^{q+1}+r-\lambda a)|\lambda \in\GF(q^2)\}\cup
\{(0,0,0,1)\}\]
and intersects the plane $\pi$ at
$\rho(P_{t,r})=(1,t,0,at^2+bt^{q+1}+r)$. We are going to show that
no $2$-secant lines of $\phi(\cH)$ pass  through $Q$. Let
$P_{t_1,r_1}(1,t_1, t_1^2, bt_1^{q+1}+r_1)$ and
$P_{t_2,r_2}(1,t_2, t_2^2, bt_2^{q+1}+r_2)$ be two distinct points
of $\phi(\cH)$. The line   $P_{t_1,r_1} P_{t_2,r_2}$ is the point
set
\[ \{ ( \lambda+1,t_1+\lambda t_2, t_1^2+ \lambda t_2^2, b(t_1^{q+1}+\lambda t_2^{q+1})+r_1+\lambda r_2)|\lambda \in\GF(q^2)\}\cup \{P_{t_2,r_2} \}. \]
If the point $Q$ were  on the line   $P_{t_1,r_1} P_{t_2,r_2}$ then
 $\lambda=-1$,  $t_1-t_2=0$ and $t_1^2-t_2^2\neq 0$, which is
impossible. Therefore, $|\rho(\phi(\cH))|=q^3+1$ and it is possible to
choose
homogeneous coordinates for the plane $\pi$ in such a way as $
\rho(\phi(\cH))$ turns out to be the set
\[\{(1,t,at^2+bt^{q+1}+r)|t\in\GF(q^2), r \in\GF(q)\}\cup \{P_\infty\};\]
that is, $\rho(\phi(\cH))$ is a non--classical B-M unital in $\pi$.

\end{document}